\documentclass[12pt]{article}
\usepackage{amsmath}
\usepackage{amssymb}
\usepackage{graphicx}
\usepackage[top=1in, bottom=0.5in, left=1in, right=1in,includefoot]{geometry}
\usepackage{enumerate}
\usepackage{algorithm,algorithmic}
\usepackage{hyperref}
\usepackage[T1]{fontenc}
\usepackage[utf8]{inputenc}
\usepackage{mathptmx}
\usepackage{mathrsfs}
\usepackage{amsthm}

\newtheorem{definition}{Definition}[section]
\newtheorem{proposition}{Proposition}[section]
\newtheorem{corollary}{Corollary}[section]
\newtheorem{remark}{Remark}[section]
\floatname{algorithm}{Procedure}

\newcommand{\algorithmicbreak}{\textbf{break}}
\newcommand{\BREAK}{\STATE \algorithmicbreak}
\usepackage{graphicx,caption,subcaption}
\usepackage{tikz}
\usetikzlibrary{intersections,arrows,automata,er,calc,backgrounds,mindmap,folding,patterns,decorations.markings,fit,snakes,shapes,matrix,positioning,calendar,topaths,trees,shapes.geometric,through}
\tikzstyle{arrow}=[thick,<-->,>=stealth]
\tikzstyle{vertex} = [circle,draw,inner sep=1pt,minimum size=4pt]
\newcommand{\vertex}{\node[vertex]}
\usepackage{pgf,pgfplots}

\begin{document}
	
	\title{Partial Domination in Prisms of Graphs}
\date{}
\author{ L. Philo Nithya$^{1}$ and Joseph Varghese Kureethara$^{1,2}$\\
	$^{1}$\small{Department of Mathematics, Christ University, Bengaluru, India.}}
	
	
	\maketitle
	
	\begin{abstract}
		For any graph G = (V, E) and proportion $p\in(0,1]$,  a set $S\subseteq V$ is a p-dominating set if $\frac{|N[S]|}{|V|}\geq p$. The $p$-domination number $\gamma_{p}(G)$ equals the minimum cardinality of a $p$-dominating set in G. For a permutation $\pi$ of the vertex set of G, the graph $\pi$G is obtained from two disjoint copies $G_1$ and $G_2$ of $G$ by joining each v in $G_1$ to $\pi(v)$ in $G_2$. i.e., $V(\pi G)= V(G_1)\cup V(G_2) \text{ and } E(G)= E(G_1)\cup E(G_2)\cup \{\{v,\pi(v)\}: v\in V(G_1), \pi(v)\in V(G_2)\}$. The graph $\pi G$ is called the prism of $G$ with respect to $\pi$. In this paper, we find some relations between the domination and the $p$-domination numbers in the context of graph and its prism graph for particular values of $p$.\\ 
		\textbf{Keywords/Phrases}:
		Permutation graph, algebraic graph theory, prism graph
		
	\end{abstract}

	\section{Introduction}
	
	The concept of prisms of graphs was first introduced by Chartrand and Harary $\cite{Chartrand1967}$ in 1967. They used the term \textit{permutation graphs} to define such graphs; but their definition was different from the one we have for permutation graphs as defined in $\cite{Lempel1972}$. Later those graphs were named as prisms of graphs with respect to a permutation. Prisms of graphs play a great role in  designing computer networks. 
	
	Partial domination \cite{Case2017, Das2018} in graphs is a variation of domination introduced in 2017. In \cite{Nithya2020}, we see some algebraic properties of the partial dominating sets of a graph. Here, in this paper we study prism graphs in the context of partial domination. 
	
	\section{Basic Terminologies}
	Let $G=(V(G),E(G))$ be a finite, simple and undirected graph with $V(G)$ as its vertex set and $E(G)$ as its edge set. For any $v \in V(G)$, $N_G(v)=\{u\in V(G):uv\in E(G)\}$ and $N_G[v]=N_G(v)\cup \{v\}$ denote the open and the closed neighborhoods of v respectively. A set $S \subseteq V(G)$ is called a dominating set if every vertex in $V-S$ is adjacent to at least one vertex in $S$. The minimum cardinality of a dominating set is called the domination number and is denoted by $\gamma(G)$. 
	
	\begin{figure*}[t!]
		\centering
		\begin{subfigure}[b]{0.5\textwidth}
			\centering
			\includegraphics[height=1.2in]{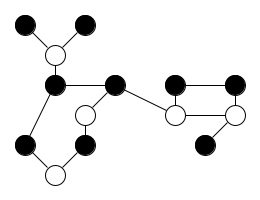}
			\caption{Domination}
		\end{subfigure}%
		~ 
		\begin{subfigure}[b]{0.5\textwidth}
			\centering
			\includegraphics[height=1.2in]{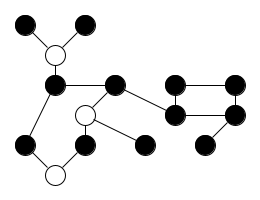}
			\caption{$\frac{2}{3}$-domination}
		\end{subfigure}
		\caption{Domination and Partial Domination}\label{1b}
	\end{figure*}

	A subset $S\subseteq V$ is called a p-dominating set for $p \in (0,1]$ if $\frac{|N[S]|}{|V|}\geq p$. The p-domination number, denoted by $\gamma_{p}(G)$ is the cardinality of the minimum p-dominating set.

	In Figure \ref{1b} (a), the white vertices dominate all the vertices of the graph. In Figure \ref{1b} (b), the white vertices do not dominate the vertices of the graph. However, they dominate exactly 10 vertices of the graph. Hence, we say, the white vertices $\frac{2}{3}$-dominate the graph.
	
	Let $\pi$ be any permutation on $V(G)$. The prism $\pi G$ of $G$ with respect to $\pi$ is obtained by taking two disjoint copies of $G$ and joining each vertex v in one copy of $G$ with $\pi(v)$ in the other copy by means of an edge.  The study of domination in prisms started in 2004 \cite{Burger2004}. Since then many works related to various domination parameters were studied \cite{Chaluvaraju2015,Hurtado2017}. 
	
	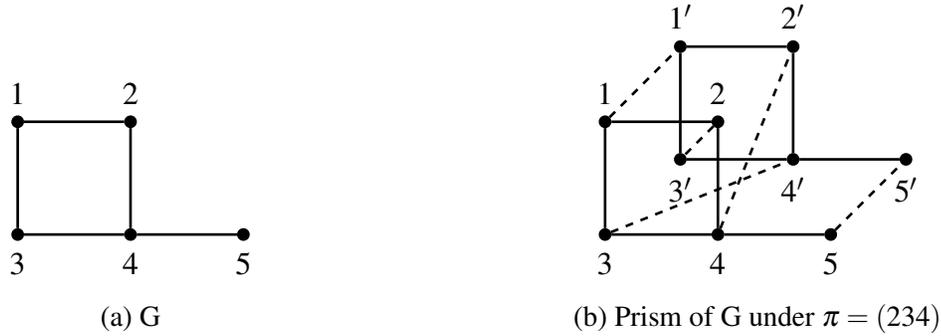
\begin{figure}
		\begin{subfigure}[b]{0.5\textwidth}
			\centering
			\begin{tikzpicture}[auto,node distance=2.5cm, thick,main node/.style={circle,draw,font=\sffamily\Large\bfseries},scale=0.5] 
				\vertex (3) at (0,0) [fill,label={below:$3$}]{}; 
				\vertex (4) at (3,0) [fill,label={below:$4$}]{}; 
				\vertex (2) at (3,3) [fill,label={above:$2$}]{}; 
				\vertex (1) at (0,3) [fill,label={above:$1$}]{}; 
				\vertex (5) at (6,0) [fill,label={below:$5$}]{}; 
				\path[every node/.style={font=\sffamily\small},line width=1pt]
				(1) edge (2) 
				(3) edge (4) 
				(4) edge (2) 
				(3) edge (1) 
				(4) edge (5) 
				; 
			\end{tikzpicture}\caption{ G}
		\end{subfigure}
		\begin{subfigure}[b]{0.5\textwidth}
			\centering
			\begin{tikzpicture}[auto,node distance=2.5cm, thick,main node/.style={circle,draw,font=\sffamily\Large\bfseries},scale=0.5] 
				\vertex (3) at (0,0) [fill,label={below:$3$}]{}; 
				\vertex (4) at (3,0) [fill,label={below:$4$}]{}; 
				\vertex (2) at (3,3) [fill,label={above:$2$}]{}; 
				\vertex (1) at (0,3) [fill,label={above:$1$}]{}; 
				\vertex (5) at (6,0) [fill,label={below:$5$}]{}; 
				\vertex (3') at (2,2) [fill,label={below:$3'$}]{}; 
				\vertex (4') at (5,2) [fill,label={below:$4'$}]{}; 
				\vertex (2') at (5,5) [fill,label={above:$2'$}]{}; 
				\vertex (1') at (2,5) [fill,label={above:$1'$}]{}; 
				\vertex (5') at (8,2) [fill,label={below:$5'$}]{}; 
				\path[every node/.style={font=\sffamily\small},line width=1pt]
				(1) edge (2) 
				(3) edge (4) 
				(4) edge (2) 
				(3) edge (1) 
				(4) edge (5)
				(1') edge (2') 
				(3') edge (4') 
				(4') edge (2') 
				(3') edge (1') 
				(4') edge (5')
				(1) edge[dashed] (1') 
				(3) edge[dashed] (4') 
				(4) edge[dashed] (2') 
				(2) edge[dashed] (3') 
				(5) edge[dashed] (5')
				; 
			\end{tikzpicture}\caption{ Prism of G under $\pi=(234)$}
		\end{subfigure}
		\caption{Prism of Graph}\label{fig:pgraph1}
	\end{figure}	
	
	In \cite{Gu2009}, it has been proved that  for any graph G, $\gamma(G)\leq \gamma(\pi G) \leq 2\gamma(G)$. Also, it has been defined in $\cite{Burger2004}$ that,  G is called universal $\gamma$-fixer, if $\gamma(G)=\gamma(\pi G)$ for all permutations $\pi$ of V(G) and $G$ is called universal doubler if $\gamma(\pi G) = 2\gamma(G)$ for all permutations $\pi$ of V(G). Analogous to this, in the context of partial domination we give the following definition: 
	\begin{definition}
		Let $p\in [0,1]$. G is called universal $\gamma_p$-fixer, if $\gamma_p(G)=\gamma_p(\pi G)$ for all permutations $\pi$ of V(G) and $G$ is called universal $\gamma _p$-doubler if $\gamma_p(\pi G) = 2\gamma_p(G)$ for all permutations $\pi$ of V(G).	
	\end{definition}
	
	Figure \ref{fig:pgraph1} is an example of a graph $G$ and its prism with respect to the permutation $\pi$=(234).
	
	\section{ Results}  
	\begin{proposition}
		Let $G$ be any n-vertex connected graph without isolated vertices and with $\gamma$=1. Then for any permutation $\pi $ of V(G) and for any $p\in (0,1]$,\\
		$\gamma_p(\pi G)$= 
		$	\begin{cases}
			1, &\text{for } p\in (0,\frac{n+1}{2n}]\\
			2, &\text{for } p\in (\frac{n+1}{2n},1]
		\end{cases}$
		
	\end{proposition}
	\begin{proof}
		Let $\pi$ be any permutation of V(G). Consider $\pi$G.\\
		Since $\gamma$(G)=1, $\exists v\in V(G)	$ such that deg(v)=n-1. This v dominates n+1 vertices in $\pi$G.\\
		Hence $\frac{|N_{\pi G}[v]|}{|V(\pi G)|}=\frac{n+1}{2n}$.\\
		Thus $\gamma_p(\pi G)= 1, \text{for } p\in (0,\frac{n+1}{2n}]$\\
		Now consider $S=\{v,\pi(v)\}$ in $\pi G$.\\
		This S dominates $\pi G$ and is minimum.\\
		Thus $\gamma_p(\pi G)= 2, \text{for } p\in (\frac{n+1}{2n},1]$.
	\end{proof}
	\begin{corollary}
		Let $G$ be any n-vertex connected graph without isolated vertices and with $\gamma$=1. Then $G$ is a universal $\gamma_p$-fixer for $p\in (0,\frac{n+1}{2n}]$ and is a universal $\gamma_p$-doubler for $p\in (\frac{n+1}{2n},1]$.
	\end{corollary}
	
	\begin{proposition}\label{d}
		Let $G$ be any n-vertex graph. Then for any permutation $\pi $ of V(G) and for any $p\in (0,\frac{n+\gamma(G)}{2n}]$, 
		$\gamma_p(\pi G) \leq \gamma(G)$.
	\end{proposition}
	\begin{proof}
		Let $\pi$ be any permutation of V(G). Consider $\pi$G.\\
		Let S be a $\gamma$-set of G. Then by the definition of $\pi G$, S is a $n+\gamma(G)$ dominating set in $\pi G$.\\
		Also if $p\leq q$ then $\gamma_p \leq \gamma_q$.\\
		Hence for any $p\in (0,\frac{n+\gamma(G)}{2n}]$, 
		$\gamma_p(\pi G) \leq \gamma(G)$.
		
	\end{proof}

	The following result shows that the above bound is sharp.
	
	\begin{proposition}
		Let $G$ be $P_n$ or $C_n$ for $n\geq 2$. Then for any permutation $\pi$ on V(G), $\gamma_\frac{n+\gamma(G)}{2n}(\pi G)=\gamma(G)$.
	\end{proposition}
	
	\begin{proof}
		Let $\pi$ be any permutation on V(G). Then by the previous proposition $\gamma_\frac{n+\gamma(G)}{2n}(\pi G)\leq\gamma(G)$. Hence it is enough if we prove that $\gamma_\frac{n+\gamma(G)}{2n}(\pi G)\geq\gamma(G)$ for any permutation $\pi$ on V(G).
		Let us assume the contradiction that $\gamma_\frac{n+\gamma(G)}{2n}(\pi G)<\gamma(G)$ for some permutation $\pi$.\\
		WLG let $S\subseteq V(G)$ be a $\gamma_\frac{n+\gamma(G)}{2n}(\pi G)$-set with $\gamma(G)-1$ vertices. Let $G_1$ and $G_2$ denote the two copies of $G$ in $\pi G$. Then two cases may arise.\\
		\textbf{Case(i)}: All the vertices of S are from either $V(G_1)$ or $V(G_2)$\\
		In this case S dominates atmost $(n-1)+(\gamma(G)-1 )$ vertices in $\pi G$.\\
		\begin{flalign*}
			\implies\frac{|N[S]|}{2n}&\leq \frac{n-1+\gamma(G)-1}{2n}&\\
			&=\frac{n+\gamma(G)-2}{2n}&\\
			&<\frac{n+\gamma(G)}{2n}&
		\end{flalign*}
		This is a contradiction to our assumption that $S\subseteq V(G)$ is a $\gamma_\frac{n+\gamma(G)}{2n}(\pi G)$-set. Hence the proof in this case.\\
		\textbf{Case(ii)}: Vertices of S are from both $G_1$ and $G_2$.\\
		WLG let there be l and m vertices from $G_1$ and $G_2$ respectively, where $l+m=\gamma(G)-1$ by our assumption. Then we have the following:\\
		
		\begin{flalign}
			\frac{|N[S]|}{2n}&\leq \frac{4(l+m)}{2n}&\nonumber\\
			&=\frac{4(\gamma(G)-1)}{2n}&\label{s}\\
			\text{ Now, }\gamma(G)&=\left\lceil\frac{n}{3}\right\rceil&\nonumber\\
			&<\frac{n}{3}+1&\nonumber\\
			\implies 3\gamma(G)&<n+3&\nonumber\\
			&<n+4&\label{x}
		\end{flalign}
		
		Hence from equations (\ref{s}) and (\ref{x}) we will get a contradiction to our assumption. Hence the proof.
	\end{proof}
	\begin{remark}
		For any graph G, $\gamma_\frac{n+\gamma(G)}{2n}(\pi G)=1$ if and only if $\gamma(G)=1$.
	\end{remark}
	
	\begin{proposition}
		Let $G$ be any graph with $u_\Delta$ as a vertex having the maximum 
		degree $\Delta(G)$ and $u'_\Delta$ as its mirror image in the second copy of $G$ in $\pi G$. Then $\gamma_\frac{n+\gamma(G)}{2n}(\pi G)=2$ if and only if one of the following two conditions holds for G:\\
		(i) $\gamma(G)=2$\\
		(ii)$\gamma(G)\geq 3$ and $\Delta(G)\geq \frac{n+\gamma(G)-4+i}{2}$ where  $|N[u_\Delta]\cap N[u'_\Delta]|=i$ for $0\leq i \leq 2$.
	\end{proposition}
	\begin{proof}
		Let us assume that $\gamma_\frac{n+\gamma(G)}{2n}(\pi G)=2$. Then by the above remark $\gamma(G)\geq 2$. If $\gamma(G)=2$,then the result is true. Hence we assume that $\gamma(G)\geq 3$. Let $|N[u_\Delta]\cap N[u'_\Delta]|=i$ for $0\leq i \leq 2$.\\
		We prove by the method of contradiction. Suppose $\Delta(G)< \frac{n+\gamma(G)-4+i}{2}$.\\
		Then $|N_{\pi G}[u_\Delta]\cup N_{\pi G}[u'_\Delta]|<n+\gamma(G)$ which implies that $\gamma_\frac{n+\gamma(G)}{2n}(\pi G)>2$ which is a contradiction. Hence the condition is necessary.\\	
		For the proof of the sufficient part, let us assume that $\gamma(G)=2$. Let S be a $\gamma(G)$-set of G. Then S is also a $\gamma_\frac{n+\gamma(G)}{2n}(\pi G)$-set of $\pi G$. Hence $\gamma_\frac{n+\gamma(G)}{2n}(\pi G)=2$.\\
		Now, let us assume the condition(ii).\\
		Consider $T=\{u_\Delta,u'_\Delta\}$.Then 
		\begin{flalign*}
			N_{\pi G}[T]&\geq 2\left[\frac{n+\gamma(G)-4+i+4-i}{2}\right]&\\
			&=n+\gamma(G)&
		\end{flalign*} 
		Also T is minimum in this context. Hence $\gamma_\frac{n+\gamma(G)}{2n}(\pi G)=2$.
	\end{proof}

	\begin{proposition}
		Suppose $p\in (0,1]$ and $\pi$ is any permutation on V(G), where $G$ is any n- vertex graph. Then $\gamma_p(G)\leq \gamma_p(\pi G)\leq 2\gamma_p(G)$.
	\end{proposition}
	\begin{proof}
		Let us first prove the lower bound part.\\
		Let $p\in (0,1]$ and S be a $\gamma_p$-set in $\pi G$.
		\begin{flalign}
			\implies &\frac{|N[S]|}{2n}\geq p&\nonumber\\
			\implies &|N[S]|\geq 2np& \label{a}
		\end{flalign}
		Let $G_1$ and $G_2$ denote the two copies of $G$ in $\pi G$.\\
		Now two cases may arise.\\
		\textbf{Case(i)}: All the vertices of S are from either $V(G_1)$ or $V(G_2)$\\
		\begin{flalign*}
			\text{Then from equation (\ref{a}), }|N[S]|&\geq 2np&\\
			&>np&
		\end{flalign*}
		Hence S is a p-dominating set in $G_1$ or $G_2$.\\
		Thus $\gamma_p(G)\leq \gamma_p(\pi G)$ in this case.\\
		\textbf{Case(ii)}: Vertices of S are from both $G_1$ and $G_2$.\\
		Let $S=X \cup Y$ where $X\subseteq V(G_1)$ and $Y\subseteq V(G_2)$.\\
		By our assumption, 
		\begin{flalign}
			\frac{|N[X\cup Y]|}{2n}& \geq p& \nonumber\\
			\implies |N[X]|+|N[Y]|&\geq 2np&\label{z}	
		\end{flalign}
		Now two cases may arise:\\
		\textbf{Subcase(i)} $|N[X]|\geq np$ and $|N[Y]|\geq np$.\\
		Let $X^*=\{w\in V(G_2)/\pi(v)=w\forall v\in X\}$ and $Y^*=\{w\in V(G_1)/\pi(w)=v\forall v\in Y\}$. Then $X\cup Y^*$ and $X^*\cup Y$ are p-dominating sets of G.\\
		Thus, 
		\begin{flalign*}
			\gamma_p(G)	&\leq |X\cup Y^*|&\\
			&=|X\cup Y|&\\
			&=\gamma_p(\pi G)&
		\end{flalign*}
		\textbf{Subcase(ii)} WLG let $|N[X]|<np$. Then $|N[Y]|\geq np+(np-|N[X]|)$ by our assumption found in equation (\ref{z}).\\
		In this case $|N[Y]|>np$.\\
		Hence $X^*\cup Y$ is a p-dominating set of G, where $X^*$ is defined as in the above subcase. Thus
		\begin{flalign*}
			\gamma_p(G)	&\leq |X^*\cup Y|&\\
			&=|X\cup Y|&\\
			&=\gamma_p(\pi G)&
		\end{flalign*}
		Hence $\gamma_p(G)\leq \gamma_p(\pi G)$ in both the cases.\\
		Now let us prove the upper bound part.\\
		Let S be a $\gamma_p$ set of G. Let $S^*=\{v^*\in G_2/\pi (v)=v^*\forall v\in S\}$.\\
		Then $\frac{|N[S\cup S^*]|}{2n}\geq p$.\\
		Thus $S\cup S^*$ is a p-dominating set of $\pi G$ and $|S\cup S^*|$ =$2\gamma_p(G)$.\\
		Hence $\gamma_p(\pi G)\leq 2\gamma_p(G)$.
	\end{proof}
	
	\begin{proposition}
		Let $G$ be a graph having an independent set $M=\{v_1, v_2,..., v_k\}$ of k-vertices $(k \geq 1)$ each having the maximum degree $\Delta(G)$. If  $N(v_i)\cap N(v_j)=\phi \forall v_i, v_j\in M$ then for $i=1,2,...,k$, $\gamma_p(\pi G)=i$ for $p\in \left (\frac{(i-1)(\Delta(G)+2)}{2n},\frac{i(\Delta(G)+2)}{2n}\right ]$ for any permutation $\pi$ on V(G).
	\end{proposition}
	
	\begin{proof}
		
		Let $N(v_i)\cap N(v_j)=\phi \forall v_i, v_j\in M$.\\
		In this case, for $1\leq i \leq k$, $S=\{v_1,v_2,...,v_i\}\subseteq M$ dominates $i(\Delta(G)+2)$ vertices in $\pi G$. Also since each vertex in M is of maximum degree, `i' is the minimum number of vertices that are required to dominate $i(\Delta(G)+2)$ vertices in $\pi G$.\\
		Hence $\gamma_\frac{i(\Delta(G)+2)}{2n}=i$.\\
		By the same argument we can say that $\gamma_\frac{(i-1)(\Delta(G)+2)}{2n}=i-1$.\\
		Also, $\frac{(i-1)(\Delta(G)+2)}{2n}$ is the best proportion of domination possible with i-1 vertices. Hence the result is true in this case.	
		
	\end{proof}
	
	\begin{proposition}
		Let $G$ be a n-vertex graph having a set of k-mutually non-adjacent vertices$(k \geq 1)$ say $M=\{v_1, v_2,..., v_k\}$ each having the maximum degree $\Delta(G)$. Let $G$ and $G'$ be the two copies of $G$ in $\pi G$ for any permutation $\pi$ on V(G). Let $M'=\{v_1', v_2',..., v_k'\}$ be the copy of M in $G'$. If for each $v_r\in M$, there exists exactly one $v_s\in M$ such that $|N(v_r)\cap N(v_s)|=1$, then for $i=1,2,...,k$, $\gamma_p(\pi G)=i$ for $p\in \left (\frac{(i-1)(\Delta(G)+2)}{2n},\frac{i(\Delta(G)+2)}{2n}\right ]$ for the following permutations $\pi$ on V(G):\\	
		\textbf{(i)} $\pi$=1\\
		\textbf{(ii)} $\pi(v_r)=v_s'$ and $\pi(v_s)=v_r'$ for each $v_r$ and $v_s$ as defined above.\\
		\textbf{(iii)} $\pi \neq 1$ and $\pi(v_r)\neq v_s'$ or $\pi(v_s) \neq v_r'$ or both and $\pi(v_i)\in M'$ $ \forall v_i \in M$.	
	\end{proposition}
	\begin{proof}
		We shall prove the theorem in three cases. For $1\leq i \leq k$, WLG assume that i is an even number. Let $S=\{v_1,v_2,...,v_i\}\subseteq M$ be such that there exists $\frac{i}{2}$pairs of $v_r,v_s$ as defined above.\\	
		\textbf{Case(i):}  $\pi=1$\\
		In this case $\{v_r,v_s'\}$ for each pair of $v_r,v_s$ will  dominate $\frac{2(\Delta(G)+2)}{2n}$ vertices in $\pi G$. Thus there exists $\frac{i}{2}$ such pairs which dominate $\frac{i(\Delta(G)+2)}{2n}$ vertices in $\pi G$. Hence by an argument similar to the previous proposition, the result is proved in this case.\\
		\textbf{Case(ii):} $\pi(v_r)=v_s'$ and $\pi(v_s)=v_r'$ for each $v_r$ and $v_s$ as defined above.\\
		In this case $\{v_r,v_r'\}$ or $\{v_s,v_s'\}$ for each pair of $v_r,v_s$ will serve the purpose for the required result.\\
		\textbf{Case(iii):} $\pi \neq 1$ and $\pi(v_r) \neq v_s'$ or $\pi(v_s) \neq v_r'$ or both and $\pi(v_i) \in M'$ $\forall v_i\in M$\\
		For this case we give an algorithm which returns a set $T$ of k-vertices from $M\cup M'$ whose members are mutually non-adjacent to each other.\\
		Let $\pi(v_i)=v_i^*$ $\forall v_i \in M$ under the above permutations and $m=k$.	
		\begin{algorithm}[H]
			\caption{Algorithm to find T}
			\begin{algorithmic}[1]
				\REQUIRE $m$,$\forall i=1$ to $m$ $v_i,{v_i}^*,N(v_i),N({v_i}^*)$ 
				\ENSURE T
				\STATE {T=$\{\}$, $i=1$,$j=1$,$k=1$}
				\WHILE {$i \leq m$}
				\IF {$\{v_i,v_i^*\}\cap T=\phi$}
				\STATE $T=T\cup \{v_i\}$
				\BREAK
				\ELSE 
				\IF {$|T|=m$}
				\RETURN T
				\ELSE 
				\STATE $i=i+1$ and go to 2
				\ENDIF 
				\ENDIF
				\ENDWHILE
				\FORALL {$j \neq i$ and $j=1$ to $m$}
				\IF {$|N(v_i) \cap N(v_j)|=1$ }
				\STATE $T=T \cup {v_j^*}$
				\BREAK
				\ENDIF
				\ENDFOR
				\FORALL {$k \neq j$ and $k=1$ to $m$}
				\IF {$|N(v_j^*) \cap N(v_k^*)|=1$}
				\STATE $i=k$ and go to 2			
				\ENDIF
				\ENDFOR
				
			\end{algorithmic}		
		\end{algorithm}	
		
		Each vertex in $T$ thus got from the algorithm  will  dominate $\frac{\Delta(G)+2}{2n}$ vertices in $\pi G$. And for $1 \leq i \leq k$, $i$ vertices from $T$ will dominate $\frac{i(\Delta(G)+2)}{2n}$ vertices in $\pi G$.		
	\end{proof}
	
	\bibliographystyle{IEEEtran}
	\bibliography{prism-graph}

\begin{thebibliography}{1}
\providecommand{\url}[1]{#1}
\csname url@samestyle\endcsname
\providecommand{\newblock}{\relax}
\providecommand{\bibinfo}[2]{#2}
\providecommand{\BIBentrySTDinterwordspacing}{\spaceskip=0pt\relax}
\providecommand{\BIBentryALTinterwordstretchfactor}{4}
\providecommand{\BIBentryALTinterwordspacing}{\spaceskip=\fontdimen2\font plus
\BIBentryALTinterwordstretchfactor\fontdimen3\font minus
  \fontdimen4\font\relax}
\providecommand{\BIBforeignlanguage}[2]{{%
\expandafter\ifx\csname l@#1\endcsname\relax
\typeout{** WARNING: IEEEtran.bst: No hyphenation pattern has been}%
\typeout{** loaded for the language `#1'. Using the pattern for}%
\typeout{** the default language instead.}%
\else
\language=\csname l@#1\endcsname
\fi
#2}}
\providecommand{\BIBdecl}{\relax}
\BIBdecl

\bibitem{Chartrand1967}
G.~Chartrand and F.~Harary, ``{Planar Permutation Graphs},'' \emph{Ann. Inst.
  Henri Poincar{\'{e}}}, vol.~4, pp. 433--438, 1967.

\bibitem{Lempel1972}
S.~Even, A.~Pnueli, and A.~Lempel, ``{Permutation Graphs and Transitive
  Graphs},'' \emph{Journal of the Association for Computing Machinery},
  vol.~19, pp. 400--410, 1972.

\bibitem{Case2017}
B.~Case, S.~Hedetniemi, R.~Laskar, and D.~Lipman, ``{Partial Domination in
  Graphs},'' \emph{Congressus Numerantium}, vol. 228, pp. 85--95, 2017.

\bibitem{Das2018}
A.~Das, ``{Partial Domination in Graphs},'' \emph{Iranian Journal of Science
  and Technology, Transactions A: Science}, pp. 1--6, 07 2018.

\bibitem{Nithya2020}
L.~P. Nithya and J.~V. Kureethara, ``{On Some Properties of Partial Dominating
  Sets},'' \emph{AIP Conference Proceedings}, 2020.

\bibitem{Burger2004}
A.~P. Burger and C.~M. Mynhardt, ``{On the Domination Number of Prisms of
  Graphs},'' \emph{Discussiones Mathematicae}, vol.~24, pp. 303--318, 2004.

\bibitem{Chaluvaraju2015}
B.~Chaluvaraju and C.~Appajigowda, ``{The Split Domination Number of a Prism
  Graph},'' \emph{Advances and Applications in Discrete Mathematics}, vol.~16,
  pp. 67--76, 2015.

\bibitem{Hurtado2017}
F.~Hurtado and M.~Mora, ``{Distance 2-Domination in Prisms of Graphs},''
  \emph{Discussiones Mathematicae}, vol.~37, pp. 383--397, 2017.

\bibitem{Gu2009}
W.~Gu and K.~Wash, ``{Bounds on the Domination Number of Permutation Graphs},''
  \emph{Journal of Interconnection Networks}, vol.~10, no.~3, pp. 205--217,
  2009.

\end{thebibliography}

\end{document}